\documentclass[a4paper,12pt]{article}
\usepackage{amsmath}
\usepackage{amsthm}
\usepackage{amssymb}
\usepackage{amscd}
 \makeatletter

 \@addtoreset{equation}{section}
  \makeatother

\theoremstyle{plain}
\newtheorem{thm}{Theorem}[section]
\newtheorem{prop}[thm]{Proposition}
\newtheorem{cor}[thm]{Corollary}
\newtheorem{lem}[thm]{Lemma}
\theoremstyle{definition}
\newtheorem{defn}{Definition}[section]

\theoremstyle{remark}
\newtheorem{rem}{Remark}[section]

\begin{document}
\title{The electromagnetic aspect for Yang-Mills fields
}
\author{ Tosiaki Kori\\
Department of Mathematics\\
School of Science and Engineering\\
 Waseda University \\
3-4-1 Okubo, Shinjuku-ku
Tokyo, Japan.\\ e-mail: kori@waseda.jp
}

\date{ }
\maketitle

\begin{abstract} 
 Let \(\mathcal{A}\) be the space of irreducible connections ( vector potentials ) over the principal bundle \(M\times SU(n)\) on a compact three-dimensional manifold \(M\).   
The Yang-Mills field  \(\mathcal{F}\) is defined as a subspace of the Whitney's direct sum  \(\mathbb{T}=T\mathcal{A}\times_{\mathcal{A}}T^{\ast}\mathcal{A}\) of the tangent and cotangent bundles of \(\mathcal{A}\).    We shall prove the Maxwell equations
\begin{eqnarray}
d_A^{\ast}B\,+\,\dot E\,=\,0\,&,&\quad d_AE\,-\,\dot B\,=\,0,\\[0.2cm]
d_AB\,=\,0\,\,&,&\quad d_A^{\ast}E\,=0\,,
\end{eqnarray}
on \(\mathcal{F}\).   Where the point of \(\mathbb{T}\) is enoted by \((A,E,B)\) with \(E\in T_A\mathcal{A}\simeq \Omega^1(M,su(n))\) and \(B\in T^{\ast}_A\mathcal{A}\simeq \Omega^2(M,su(n))\).    The first two equations are the  Hamilton equations of motion derived from a symplectic structure on \(\mathbb{T}\),  and  the second equations which are the defining equations of \(\mathcal{F}\) come from the action of the group of gauge transformations \(\mathcal{G}\) on \(\mathcal{A}\).    
The symplectic structure on \(\mathbb{T}\) is given by the 2-form:
 \begin{eqnarray*}
\Omega_{(A,E,B)}\left(\left(\begin{array}{c}a_1\\e_1\\ \beta_1\end{array}\right)\,,\,\left(\begin{array}{c}a_2\\e_2\\ \beta_2\end{array}\right)\,\right)&=&\,(e_2\wedge a_1\,,B\,)\,-\,(e_1\wedge a_2\,,B)\,\\[0.2cm]
\qquad&\,&\,+
\left(\,e_2\,,\,d_A^{\ast}\beta_1\,\right)_1\,-\,\left(\,e_1\,,\,d_A^{\ast}\beta_2\,\right)_1\,,
\end{eqnarray*}
for \((A,E,B)\in \mathbb{T}\) and \(\,\left(\begin{array}{c}a_i\\e_i\\ \beta_i\end{array}\right) \in T_{(A,E,B)}\mathbb{T}\,\), \(i=1,2\), with \(a_i,\,e_i\in T_A\mathcal{A}\) and \(\beta_i\in T_A^{\ast}\mathcal{A}\).   The corresponding Poisson bracket on \(\mathcal{F}\) is  
\begin{equation}
\left\{\,\Phi\,,\,\Psi\,\right\}_{(E,B)}^{\mathbb{T}}\,=\,\left(\,\frac{\delta\Phi}{\delta B}\,,\,d_A^{\ast}\frac{\delta\Psi}{\delta E}\,\right)_1\,-
\,\left(\,\frac{\delta\Psi}{\delta B}\,,\,d_A^{\ast}\frac{\delta\Phi}{\delta E}\,\right)_1\,.
\end{equation}
This is a parallel formula of Marsden-Weinstein  in case of the electric-magnetic field.    We shall investigate the Clebsch parametrization of the Yang-Mills field \(\,(\mathcal{F},\,\Omega)\).   
 We show that the action of   \(\,\mathcal{G}\)  on \(\,(\mathcal{F},\Omega )\) is Hamiltonian with the moment map
\(
\mathbb{J}(\,E,\,B\,)\,=\,[\,d_A\ast B\,,\,E\,]\,
\).    This gives a conserved quantity  \(\int_M\,[\,d_A\ast B\,,\,E\,]\,\) which is due to the non-commutativity of the gauge group.
\end{abstract}


Subj. Class: Global analysis, Quantum field theory.


\section{Introduction}

We shall investigate the electric-magnetic paradigm on the Yang-Mills field \(\mathcal{F}\) over a compact three-dimensional manifold \(M\).    We shall give the Maxwell equations:
\begin{eqnarray}
d_A^{\ast}B\,+\,\dot E\,=\,0\,&,&\quad d_AE\,-\,\dot B\,=\,0,\label{1st}\\[0.2cm]
d_AB\,=\,0\,\,&,&\quad d_A^{\ast}E\,=0\,,\label{2nd}
\end{eqnarray}
at \((A,E,B)\in\mathcal{F}\).
The first equations (\ref{1st}) are the  Hamilton equations of motion derived from a symplectic structure on \(\mathcal{F}\),  and  the second equations  (\ref{2nd})  that represent the conservation of  electric and magnetic charges come from the action of the group of gauge transformations.    

 Let \(\mathcal{A}\) be the space of irreducible connections ( vector potentials ) on \(M\times SU(n)\) and let 
\(T\mathcal{A}\simeq \Omega^1(M,su(n))\) be the tangent space and \(T^{\ast}\mathcal{A}\simeq \Omega^2(M,su(n))\) be the cotangent space over \(\mathcal{A}\).     
Our Yang-Mills field \(\mathcal{F}\) is realized as a subspace of  the Whittney's direct sum 
 \[\mathbb{T}=T\mathcal{A}\times_{\mathcal{A}}T^{\ast}\mathcal{A}\longrightarrow \mathcal{A}.\]
       A points of \(\mathbb{T}\) is denoted by \((A,E,B)\) with \(E\in T_A\mathcal{A}\) and \(B\in T^{\ast}_A\mathcal{A}\).    Then the subspace \(\mathcal{F}\) of \(\mathbb{T}\) is defined by the equations 
(\ref{2nd}).    

\(\mathbb{T}\) becomes a symplectic manifold endowed with the following symplectic form:
 \begin{eqnarray*}
\Omega_{(A,E,B)}\left(\left(\begin{array}{c}a_1\\e_1\\ \beta_1\end{array}\right)\,,\,\left(\begin{array}{c}a_2\\e_2\\ \beta_2\end{array}\right)\,\right)&=&\,(e_2\wedge a_1\,,B\,)_2\,-\,(e_1\wedge a_2\,,B)_2\,\\[0.2cm]
\qquad&\,&\,+
\left(\,e_2\,,\,d_A^{\ast}\beta_1\,\right)_1\,-\,\left(\,e_1\,,\,d_A^{\ast}\beta_2\,\right)_1\,,
\end{eqnarray*}
for \(\,\left(\begin{array}{c}a_i\\e_i\\ \beta_i\end{array}\right) \in T_{(A,E,B)}\mathbb{T}\,\), \(i=1,2\), where 
the bracket \((\,,\,)_k\) is the inner product on the Sobolev space of differential \(k\)-forms; \(\Omega^k_s(M,\,su(n))\).    \(\mathcal{F}\) is a symplectic subspace.   
The corresponding Poisson bracket on \(\mathbb{T}\), so on \(\mathcal{F}\), becomes 
\begin{equation}
\left\{\,\Phi\,,\,\Psi\,\right\}_{(E,B)}^{\mathbb{T}}\,=\,\left(\,\frac{\delta\Phi}{\delta B}\,,\,d_A^{\ast}\frac{\delta\Psi}{\delta E}\,\right)_1\,-
\,\left(\,\frac{\delta\Psi}{\delta B}\,,\,d_A^{\ast}\frac{\delta\Phi}{\delta E}\,\right)_1\,.
\end{equation}
This is the Yang-Mills counterpart to the Poisson bracket on the electro-magnetic field \(\mathcal{F}_{Max}\) discussed by Marsden-Weinstein, \cite{A, MW2}:
\begin{equation}
\left\{\,\Phi\,,\,\Psi\,\right\}_{(E,B)}\,=\,\left(\,\frac{\delta\Phi}{\delta B}\,,\,curl\,\frac{\delta\Psi}{\delta E}\,\right)_1\,-
\,\left(\,\frac{\delta\Psi}{\delta B}\,,\,curl\,\frac{\delta\Phi}{\delta E}\,\right)_1\,.
\end{equation}
As is well known the Maxwell's equation is given by 
\begin{eqnarray}
curl\,B\,+\,\dot E\,=\,0\,&,&\quad curl\,E\,-\,\dot B\,=\,0,\\[0.2cm]
div\,B\,=\,0\,\,&,&\quad div\,E\,=0\,.\label{charge}
\end{eqnarray}
The electric-magnetic field \(\mathcal{F}_{Max}\) is the subspace characterized by the second line equations (\ref{charge}).

 We shall give a detailed explanation of our investigation.     Let \(M\) be a compact \(m\)-dimensional manifold.   
 Let \(\mathcal{A}\) be the space of irreducible connections ( vector potentials ) over the trivial principal bundle \(M\times G\) and let 
\(T\mathcal{A}\) be the tangent space and \(T^{\ast}\mathcal{A}\) be the cotangent space over \(\mathcal{A}\). 
  \(T_A\mathcal{A}\) is a vector space isomorphic to \(\Omega^1_s(M, Lie\,G) \); the Sobolev space of differential \(1\)-forms on \(M\), and \(T^{\ast}_A\mathcal{A}\) is a vector space isomorphic to \(\Omega^{m-1}_s(M, Lie\,G) \).       The dual coupling is given by 
\[\int_M\,tr\,a\wedge\beta\,,\quad a\in T_A\mathcal{A},\,\beta\in T^{\ast}_A\mathcal{A}.\]
In the following we shall abreviate to write \(R=T\mathcal{A}\) and \(S=T^{\ast}\mathcal{A}\).    But we retain the freedom to use the original notations.    A point of \(S\) is denoted by \((A,\lambda)\) with \(A\in\mathcal{A}\) and \(\lambda\in T^{\ast}_A\mathcal{A}\).   A tangent vector to  \(S\) at \((A,\lambda)\) is \(\left(\begin{array}{c}a\\ \alpha\end{array}\right)\) with \(a\in T_A\mathcal{A}\) and \(\alpha\in T^{\ast}_A\mathcal{A}\).   
   There exists a canonical \(1\)-form \(\theta\) on \(S\) characterized by \(\varphi^{\ast}\theta=\varphi\) for any section \(\varphi\) of 
\(\,T^{\ast}\mathcal{A}\longrightarrow \mathcal{A}\).    The exterior derivative ( on \(\mathcal{A}\) ) of \(\theta\) gives a canonical symplectic form \(\omega=\tilde d\theta\) on \(S\).     Then the Hamiltonian vector vector field of a function  \(\Phi=\Phi(A,\lambda)\)    
is given by    
\begin{equation*}
X_{\Phi}\,=\,\left(\begin{array}{c}\, \,\frac{\delta\Phi}{\delta \lambda}\\[0.2cm]
-\, \frac{\delta\Phi}{\delta A}\end{array}\right).
\end{equation*}
 Where \(\frac{\delta\,}{\delta A}\)  indicates the partial derivative ( in the sense of Frechet-Gateau ) to the direction of \(\mathcal{A}\), and \(\frac{\delta\,}{\delta \lambda}\) is the partial derivative along the fiber \(T^{\ast}_A\mathcal{A}\).

For example, if \(\dim M=3\) and 
\[
H(A,B)\,=\, \frac12\int_M\,tr\,F_A\wedge\ast F_A\,+\,\frac12\int_M\,tr\,B\wedge \ast B\,
\]
for \(A\in \mathcal{A}\,\), \(B\in T^{\ast}_A\mathcal{A}\), 
then the Hamiltonian vector field of \(H\) becomes 
\begin{equation*}
X_{H}\,=\,\left(\begin{array}{c}  \ast B \\[0.1cm]
-\,  d_A\ast F_A \end{array}\right),
\end{equation*}
and the Hamilton's equation of motion is  
\begin{equation}\label{*hamilton}
\dot{A}\,=\,\ast B\,,\quad \dot{B}\,=\,- \ast\,d_A^\ast F_A\,, \qquad\mbox{ for } A\in \mathcal{A}\,,\,B\in T_A\mathcal{A}\,.
\end{equation}
The group of ( pointed ) gauge transformations \(\,{\cal G}={\cal G}(M)=\Omega^0_s(M,\,Ad\,G)\) acts on \(\mathcal{A}\) ( from the right ) by 
\(\,g\cdot A=g^{-1}Ag+g^{-1}dg\,\).      \({\cal G}\) acts on 
\(T_A{\cal A}\)  by the adjoint representation; 
\(a\longrightarrow Ad_{g^{-1}}\,a=g^{-1}ag\), and on  \(T^{\ast}_A{\cal A}\) by its dual \(\alpha\longrightarrow g^{-1}\alpha g\).     The dual space of \(\,Lie\,\mathcal{G}\) is given by \( (Lie\,\mathcal{G})^{\ast}\,=\,\Omega^{m}(M, Lie\,G)\) with the dual pairing:
\begin{equation}
\langle \mu,\,\xi\rangle\,=\,\int_M\,tr\,(\mu\, \xi\,)\,,\qquad \forall\xi\in \mathcal{G}\,,\,\mu\in \Omega^m(M,Lie\,G).
\end{equation}
We see that the action of the group of gauge transformations \(\mathcal {G}\) on the symplectic space \((\,S=T^{\ast}\mathcal{A},\,\omega\,)\) is an hamiltonian action with the moment  map \(J:\,S\longrightarrow\,(Lie\,\mathcal{G})^{\ast}=\Omega^m(M,Lie\,G)\)  given by 
\begin{equation}
J^{\ast}(A,\lambda)=\,-\,d_A\,\lambda\,. 
\end{equation} 
From Marsden-Weinstein reduction theorem the reduced space \((J^{\ast})^{-1}(\rho)/\mathcal{G}\) for a \(\rho\in (Lie\,\mathcal{G})^{\ast}\) becomes a symplectic manifold endowed with the induced symplectic form \(\omega\),  and coincides with the space 
\(\{(A,\lambda)\in S\,;\, \, - d_A\lambda=\rho\,\}\).

 A parallel argument is valid on the tangent space \(R=T\mathcal{A}\) as we shall show in the following, but this is not a canonical one.    The point of  \(R\) is denoted by \((A,p)\) with \(A\in \mathcal{A}\) and \(p\in T_A\mathcal{A}\).   
The tangent space at \((A,p)\in R\) is \(T_{(A,p)}R=T_A\mathcal{A}\oplus T_A{\mathcal A}\), so any tangent vector \(\mathbf{a}\in T_{(A,p)}R\) is of the form \(\mathbf{a}=\left(\begin{array}{c}a\\x\end{array}\right)\,\) with \(a,\,x\in T_A\mathcal{A}\).   
The symplectic structure on \(R\) is defined by the formula 
\begin{equation}
\sigma_{(A,p)}\left(\,\left(\begin{array}{c}a\\x\end{array}\right)\,,\,\left(\begin{array}{c}b\\y\end{array}\right)\,\right)\,=\, (\,b\,,\,x\,)_1\,-\,(\,a\,,\,y\,)_1\,,\end{equation}
for all  \(\left(\begin{array}{c}a\\x\end{array}\right)\,,\,\left(\begin{array}{c}b\\y\end{array}\right)\,\in T_{(A,p)}R\,\).

If \(\dim \,M=3\), the hamiltonian function 
\begin{equation}
H(A,p)\,=\,\frac12(F_A,F_A)_2\,+\,\frac12(p,p)_1\,.
\end{equation}
gives the Hamilton's equation of motion 
\begin{equation}\label{hamilton}
\dot A\,=\,-p\,,\qquad 
\dot p\,=\,d_A^{\ast}F_A \,.
\end{equation}
Actualy we find that the symplectic manifolds \((S\,,\,\omega\,)\) and \((R\,,\,\sigma\,)\) are isomorphic via the Hodge operator
 \begin{equation}
 \ast: T_A\mathcal{A}\,\simeq \Omega^1(M,su(n))\,\longleftrightarrow\,
T^{\ast}_A\mathcal{A}\,\simeq \Omega^{m-1}(M,su(n))\,.\label{*iso}
\end{equation}

When \(\dim\,M=3\) the symplectic isomorphism given by Hodge \(\ast\) changes the 
Hamilton equations of motion (\ref{*hamilton}) and (\ref{hamilton})  each other:
\begin{equation}
\{\begin{array}{c} \dot{A}\,=\,\ast B\\[0.2cm]  \dot{B}\,=\,- \ast\,d_A^\ast F_A\end{array}\quad
\stackrel{\,p=\ast B\,}\Longleftrightarrow\quad
\{\begin{array}{c} \dot A\,=\,p\\[0.2cm]
\dot p\,=\,- d_A^{\ast}F_A\end{array}\,.
\end{equation}
One of our purpose in this paper is to write down the equations of motion in the form that do not contain the potential variable \(A\), but contain only  the field variables \(E\) and \(B\).
So both (\ref{*hamilton}) and (\ref{hamilton}) are insufficient for us.

   The 
action of  \(\,\mathcal{G}\) on the symplectic space \((\,R=T\mathcal{A},\sigma\,)\) is an hamiltonian action with the moment map given by
\begin{equation}
J(A,p)\,=\,d_A^{\ast}\,p\,. \label{tgmoment}
\end{equation}
This time the image of the moment map is  \((Lie\,\mathcal{G})^{\ast}\simeq Lie\,\mathcal{G}= \Omega^0(M,Lie\,G)\) by virtue of the inner product \((\, ,\,)_0\).
The symplectic reduction 
\(\left( J^{-1}(0)/\mathcal{G},\,\sigma\right)\,\) becomes a symplectic manifold and coincides with the subbundle 
 \begin{equation}\label{r0}
 R^0=\{(A,p)\in R \,;\, \,d^{\ast}_Ap=0\,\}.
  \end{equation}

Now we suppose that \(M\) is a compact manifold with \(\dim\,M=3\).    We shall introduce the space where the electric magnetic future of Yang-Mills theory is relevant.   
We consider the direct sum 
\begin{equation}
\begin{array}{ccc}
\mathbb{T}=\,R\,\times_{\mathcal{A}}\,S\,&\,\longrightarrow\,&\,(S\,,\omega )\,\\[0.3cm]
\pi_{\ast}\,\downarrow\quad  &\qquad&\,\pi \downarrow\, \\[0.4cm]
(\,R\,,\sigma)\,&\,\longrightarrow\,&\, \,\mathcal{A}\,\,.
\end{array}
\end{equation}
A point of \(\,\mathbb{T}\,\) will be denoted by \((A,E,B)\) with \(E\in R\) and \(B\in S\) that are over \(A\in\mathcal{A}\).   
The tangent space of \(\,\mathbb{T}\,\) at the point \((E,B)\in \mathbb{T}\,\) ( over \(A\in\mathcal{A}\))  is  \( T_{(A,E,B)}\mathbb{T}\,=T_A\mathcal{A}\oplus T^{\ast}_A\mathcal{A}\).   So a vector in \(T_{(A,E,B)}\mathbb{T}\,\) is denoted by \(\mathbf{a}=\left(\begin{array}{c}a\\e\\ \beta\end{array}\right) \), with \(a,\,e\in T_A\mathcal{A}\simeq \Omega^1(M,Lie\,G)\) and \(\beta\in T_A^{\ast}\mathcal{A}\simeq \Omega^2(M,Lie\,G)\).      
The inner product on the fiber \(T_{(E,B)}\mathbb{T}\) over \(A\in\mathcal{A}\) is given by 
\begin{equation}
\left( \left(\begin{array}{c}e_1\\ \beta_1\end{array}\right) ,\,\left(\begin{array}{c}e_2\\ \beta_2\end{array}\right) \right)_{\mathbb{T}}\,=\, (e_2,d_A^{\ast}\beta_1\,)_1\,+\,(e_1,\,d_A^{\ast}\beta_2\,)_1\,.
\end{equation}
The partial derivaive of a function \(\Phi=\Phi(E,B)\) over \(\mathbb{T}\) to the direction \(e\in T_A\mathcal{A}\) 
 is defined as  the vector  \(\frac{\delta\Phi}{\delta E}\,\in T_A^{\ast}\mathcal{A}\,\).   Respectively that  to the direction  \(\beta\in T^{\ast}_A\mathcal{A}\) is defined as  the vector  \(\frac{\delta\Phi}{\delta B}\,\in T_A\mathcal{A}\,\).    They satisfy the defining equation: 
 \begin{equation}\label{dPhi}
( \widetilde d\Phi)_{(E,B)}\left(\begin{array}{c}e\\\beta\end{array}\right)\,=\,
\left( \left(\begin{array}{c}\frac{\delta\Phi}{\delta B}\\[0.2cm] \frac{\delta\Phi}{\delta E}\end{array}\right) ,\,\left(\begin{array}{c}e\\[0.2cm] \beta\end{array}\right) \right)_{\mathbb{T}}\,
\end{equation}
Now we endow the space \(\mathbb{T}\) with a symplectic structure given by 
the  following 2-form \(\Omega\):
 \begin{eqnarray*}
\Omega_{(A,E,B)}\left(\left(\begin{array}{c}a_1\\e_1\\ \beta_1\end{array}\right)\,,\,\left(\begin{array}{c}a_2\\e_2\\ \beta_2\end{array}\right)\,\right)&=&\,(e_2\wedge a_1\,,B\,)_2\,-\,(e_1\wedge a_2\,,B)_2\,\\[0.2cm]
\qquad&\,&\,+
\left(\,e_2\,,\,d_A^{\ast}\beta_1\,\right)_1\,-\,\left(\,e_1\,,\,d_A^{\ast}\beta_2\,\right)_1\,,
\end{eqnarray*}
for \(\,\left(\begin{array}{c}a_i\\e_i\\ \beta_i\end{array}\right) \in T_{(A,E,B)}\mathbb{T}\,\), \(i=1,2\), where 
the bracket \((\,,\,)_k\) is the inner product on the Sobolev space of differential \(k\)-forms; \(\Omega^k_s(M,\,su(n))\).   \(\Omega\) is a non-degenerate skew-symmetric 2-form and we have a symplectic structure on \((\mathbb{T},\Omega)\).    
   For \(\Phi,\,\Psi\in\,C^{\infty}(\mathbb{T})\) the Poisson bracket is defined by the formula:
\begin{equation}
\left\{\,\Phi\,,\,\Psi\,\right\}_{(E,B)}^{\mathbb{T}}\,=\,\Omega_{(E,B)}\left(\,X_{\Phi}\,,\,X_{\Psi}\,\right)\,.
\end{equation}
We have the following representation of Poisson bracket:
\begin{equation}
\left\{\,\Phi\,,\,\Psi\,\right\}_{(E,B)}^{\mathbb{T}}\,=\,\left(\,\frac{\delta\Phi}{\delta B}\,,\,d_A^{\ast}\frac{\delta\Psi}{\delta E}\,\right)_1\,-
\,\left(\,\frac{\delta\Psi}{\delta B}\,,\,d_A^{\ast}\frac{\delta\Phi}{\delta E}\,\right)_1\,.
\end{equation}
 If we take the Hamiltonian function on \(\mathbb{T}\) written in the vortex representation;
 \begin{equation}
 H(E,B)\,=\frac12 \left(
 \left(\begin{array}{c}E\\ B\end{array}\right)\,,
 \left(\begin{array}{c}d_A^{\ast}B\\ d_AE\end{array}\right)
 \right)_{\mathbb{T}} =\frac12\{(d_AE,d_AE)_1+(d_A^{\ast}B,d_A^{\ast}B)_1\},
 \end{equation}
We obtain the following equation of motion on the strength field \(\mathbb{T}\): 
 \begin{equation}
\dot E\,=\, \frac{\delta H}{\delta B}=d_A^{\ast}B\,,\quad   \dot B\,=\frac{\delta H}{\delta E}=d_AE\,.
 \end{equation}
Which is nothing but the Maxwell equation (\ref{1st}) over \(\mathbb{T}\):

The {\it Yang-Mills field} is  the symplectic subspace of \(\mathbb{T}\) defined by 
\begin{equation}
\mathcal{F}\,=\,\left\{(E,B)\,\in  \mathbb{T}\,:\quad d_AB=0,\,\,d_A^{\ast}E=0 \quad\mbox{ with \(\,\pi_{\ast}(E,B)=A\)}\,\right\}\,
\end{equation}
It is  a \(\mathcal{G}\)-invariant subspace.  We shall give a symplectomorphism from the reduced space \(( R^0,\sigma )\), (\ref{r0}),  to \((\mathcal{F},\Omega)\), that is,   \(( R^0,\sigma )\) is  a symplectic variable ( Clebsch parametrization ) of the Yang-Mills field \(\mathcal{F}\).

We find that the action of  \(\,\mathcal{G}\) on \(\,\mathcal{F}\) is Hamiltonian with the moment map
\begin{equation}
\mathbb{J}(\,E,\,B\,)\,=\,[\,d_A\ast B\,,\,E\,]\,.
\end{equation}
Hence on the \(\mathcal{G}\)-orbit passing through a solution of  equation (\ref{1st}) 
we have \(\mathbb{J}(E,B)=\,[\,d_A\ast B\,,\,E\,]\).       We have an invariant 
\begin{equation}
\int_M\,[\,d_A\ast B\,,\,E\,]
\end{equation}
 which reflects the non-commutativity of the gauge group \(G=SU(n)\).

\section{Calculation on the space of connections \cite{DK, K, K2}}

Let \(M\) be a compact, connected and oriented \(m\)-dimensional riemannian manifold  possibly with boundary \(\partial M\).    Let \(P\stackrel{\pi}{\longrightarrow}M\) be a principal \(G\)-bundle, \(G=SU(N)\), \(N\geq 2\).

We write \({\cal A}={\cal A}(M)\) the space of {\it irreducible} \(L^2_{s-1}\) connections over \(P\), which differ from a smooth connection by a \(L^2_{s-1}\) section of \(T^{\ast}_M\otimes Lie\, G\), hence
the tangent space of \({\cal A}\) at \(A\in{\cal A}\) is
 \begin{equation}
 T_A{\cal A}=\Omega^1_{s-1}(M, Lie \,G)\,.
 \end{equation}
 The cotangent space of \({\cal A}\) at \(A\) is 
  \begin{equation}
 T_A^{\ast}{\cal A}=\Omega^{m-1}_{s-1}(M, Lie \,G)\,),
 \end{equation}
where the pairing \(\langle a,\,\alpha\rangle_A\,\) of \(\alpha\in T^{\ast}_A{\cal A}\,\) and \(a\in T_A{\cal A}\)  is given by the symmetric bilinear form \((X,Y)\longrightarrow tr(XY)\) of \(Lie\,G\) and the Sobolev norm \((\,,\,)_{s-1}\) on the Hilbert space \(L^2_{s-1}(M)\):
\[
\langle \phi\otimes X, \psi\otimes Y\rangle =(\phi\,,\, \psi\,)_{s-1}  \,tr(XY),\]
for \(
\psi\in \Omega^{m-1}(M),\, \phi\in\Omega^1(M),\) and \( X,Y\in Lie\,G\).   We shall write it by  
\(\langle a,\,\alpha\rangle_A=\int_M\,tr(\,a\wedge \alpha\,)\), or simply by  \(\int_M\,tr(\,a\,\alpha\,)\).

 A vector field  \(\mathbf{a}\) on \({\cal A}\) is a section of the tangent bundle;
 \(\mathbf{a}(A)\in T_A{\cal A}\), and a 1-form \(\varphi\) on \({\cal A}\) is a section of the cotangent bundle;
 \(\varphi(A)\in T_A^{\ast}{\cal A}\).      
    
For a smooth map \(F=F(A)\) on \({\cal A}\) valued in a vector space \(V\) the derivation \(\partial_A F\) is defined by the functional variation of \(A\in {\cal A}\):
\begin{eqnarray}
\partial_A F\,&:&\,T_A{\cal A}\longrightarrow\,V\,,\\[0.2cm]
(\partial_A F)a&=&\lim_{t\longrightarrow 0}\frac{1}{t}\left(\,F(A+ta)-F(A)\,\right),\quad \mbox{ for \(a\in T_A{\cal A}\)}.
\end{eqnarray}
For example,  
\[(\partial_A A)\,a=a,\]
since the derivation of an affine function is defined by its linear part.
The curvature of \(A\in {\cal A}\) is given by 
\[F_A=dA+\frac12[A\wedge A]\in \Omega^2_{s-2}(M, Lie\,G).\]  
So it holds that 
\[F_{A+a}=F_A+d_Aa+a\wedge a,\]
and we have
  \[(\partial_A F_A)a=d_Aa.\]
  The derivation of a vector field \(\mathbf{v}\) on \({\cal A}\) and that of a 1-form \(\varphi\) are defined similarly:
  \[(\partial_A\mathbf{v})a \in T_A{\cal A},\qquad(\partial_A\varphi)a \in T^{\ast}_A{\cal A},\qquad \forall a\in T_A{\cal A}.\]
   It follows that the derivation of a function   \(F=F(A)\) by a vector field \(\mathbf{v}\) is given by
  \[(\mathbf{v}F)_A=(\partial_AF)(\mathbf{v}_A) .\] 
  We have the following formulas,  \cite{BN, K}.
  \begin{eqnarray}
 [\, \mathbf{v},\,\mathbf{w}\,]_A&=&(\partial_A\mathbf{v})\mathbf{w}_A-(\partial_A\mathbf{w})\mathbf{v}_A,\label{derivation1}\\[0.2cm]
 ( \mathbf{v}\langle \varphi,\mathbf{u}\rangle )_A&=&\langle \varphi_A, (\partial_A\mathbf{u})\mathbf{v}_A\rangle +\langle (\partial_A\varphi)\,\mathbf{v}_A,\mathbf{u}_A\rangle .\label{derivation2}
  \end{eqnarray}
  
Let \(\widetilde d\)  be the exterior derivative on \({\cal A}(M)\).   For a function \(F\) on \({\cal A}(M)\), \((\widetilde dF)_A\,a=(\partial_AF)\,a\).   \\
 For a 1-form \( \Phi\) on \({\cal A}(M)\), 
\begin{eqnarray}
(\widetilde d \Phi)_A({\bf a},{\bf b})&=&(\partial_A< \Phi,{\bf b}>){\bf a}-(\partial_A< \Phi,{\bf a}>){\bf b}-< \Phi, [{\bf a},{\bf b}]>\nonumber \\[0.2cm]
&=& <(\partial_A \Phi){\bf a},  {\bf b}>-<(\partial_A \Phi){\bf b},{\bf a}>.\label{extder1}
\end{eqnarray}
This follows from (\ref{derivation1}) and (\ref{derivation2}).
Likewise, if \(\varphi\) is a 2-form on \({\cal A}(M)\)  then it holds that 
\begin{equation}\label{extder}
(\widetilde d\,\varphi)_A({\bf a},{\bf b},{\bf c})=(\partial_A\varphi({\bf b},{\bf c})){\bf a}+(\partial_A\varphi({\bf c},{\bf a})){\bf b}+(\partial_A\varphi ({\bf a},{\bf b})){\bf c}\,.\end{equation}

We write the group of \(L^2_s\)-gauge transformations by \({\cal G}^{\prime}(M)\):
\begin{equation}
{\cal G}^{\prime}(M)=\,\Omega^0_s(M, Ad\,P) .\end{equation}
Where \(Ad\,P=P\times_GG\) is the adjoint bundle associated to the principal bundle \(P\).   In this paper we shall mainly deal with the trivial principal bundle.   In this case 
\({\cal G}^{\prime}(M)=\,\Omega^0_s(M, G)\).    
\({\cal G}^{\prime}(M)\) acts on \({\cal A}(M)\) by 
\begin{equation}
g\cdot A=g^{-1}dg+g^{-1}Ag=A+g^{-1}d_Ag.
\end{equation}
By Sobolev lemma one sees that \({\cal G}^{\prime}(M)\) is a Banach Lie\,Group  and its action is a smooth map of Banach manifolds.   

In the following we choose a fixed point \(p_0\in M\) and deal with the group of gauge transformations that are identity at \(p_0\):
\[{\cal G}={\cal G}(M)
=\{g\in {\cal G}^{\prime}(M);\,g(p_0)=1\,\}.\] 
  \({\cal G}\) act freely on \({\cal A}\).    
 Let \({\cal C}(M)={\cal A}(M)/{\cal G}(M)\) be the quotient space of this action.   It is a smooth infinite dimensional manifold.     
 
 Let \({\cal G}_0(M)\) be the group of gauge transformations that are identity on the boundary of \(M\).       
When \(M\) has no boundary \({\cal G}(M)={\cal G}_0(M)\).

   We have 
    \[Lie\,({\cal G})=\Omega^0_s(M, ad\,P)\,.\]
    Where \(ad\,P=P\times_GLie\,G\) is the derived bumdle of \(Ad\,P\).    When \(P\) is trivial 
    \( Lie\,({\cal G})=\Omega^0_s(M, Lie\,G)\).
    The Lie algebra of \(\mathcal{G}_0\,\) is 
\[Lie\,\mathcal{G}_0\,=\{\xi\in Lie\,\mathcal{G};\, \xi\vert \partial M=0\,\}=\{\xi\in \Omega^0_s(M, ad\,P);\, \xi\vert \partial M=0\,\}\,.\]

The infinitesimal action of \({\cal G}\) on \({\cal A}\) is described by
\begin{equation}
\xi\cdot A=\,d_A\xi\,=\,d\xi\,+[A\wedge\,\xi]\, , \qquad \forall \xi\in Lie\,\mathcal{G},\,\forall A\in\mathcal{A}\,.
\end{equation}

The fundamental vector field on \({\cal A}\) corresponding to  
\(\xi \in Lie({\cal G})\) is given by 
\[d_A\xi\,=\,\frac{d}{dt}\lvert_{t=0}(\exp \,t\xi)\cdot A,\]
and the tangent space to the orbit at \(A\in{\cal A}\) is 
\begin{equation}
T_A({\cal G}\cdot A)=\{d_A\xi\,;\,\xi\in \Omega^0_s(M, ad\,P)\}.
\end{equation}

\section{Canonical structure on \(T^{\ast}{\cal A}\)}

On the cotangent bundle of any manifold we have the notion of canonical symplectic form, and the standard theory of Hamiltonian mechanics and its symmetry follows from it, \cite{AM}.   We apply these standard notions to our infinite dimensional manifold \({\cal A}(M)\) and write up their explicit formulas.    

\subsection{ Canonical 1-form and 2-form on \(T^{\ast}{\cal A}\)}

Let \(M\) be a manifold of \(\dim M=m\) possibly with the non-empty boundary \(\partial M\).   
Let \(T^{\ast}{\cal A}\stackrel{\pi}{\longrightarrow}{\cal A}\) be the cotangent bundle.    We denote the pairing of  \(T_{A}{\cal A}\) and \(T^{\ast}_{A}{\cal A}\) by 
\begin{equation}
\langle\,a\,,\,\alpha\,\rangle_A\,=\,\int_M\,tr\,a\wedge\alpha\,,\qquad \forall\,a\in T_A\mathcal{A},\,\alpha\in T^{\ast}_A\mathcal{A}\,.\end{equation}
In the following we shall denote the cotangent space \(T^{\ast}\mathcal{A}\) by \(S\).   
The point of \(S\) will be denoted by \((A,\lambda)\) with \(A\in \mathcal{A}\) and \(\lambda\in T^{\ast}_A\mathcal{A}\).  
The tangent space  to  the cotangent space \(S\,\) at the point \((A,\lambda)\in S\) becomes 
\begin{equation}T_{(A,\lambda)}S\,=T_A{\cal A}\oplus T^{\ast}_A{\cal A}=\Omega^1(M,\,Lie\,G)\oplus\Omega^{m-1}(M.\,Lie\,G).
\end{equation}
Any tangent vector \(\mathbf{a}\in T_{(A,\lambda)}S\) has the form \(\mathbf{a}=\left(\begin{array}{c}a\\\alpha\end{array}\right)\,\) with \(a\in T_A\mathcal{A}\) and  \(\alpha\in T^{\ast}_A\mathcal{A}\).   

The canonical 1-form on the cotangent space \(S\) is defined as follows:    
\begin{equation}
\theta_{(A,\lambda)}(\left(\begin{array}{c}a\\ \alpha\end{array}\right))=\langle\, \lambda,\pi_{\ast}\left(\begin{array}{c}a\\ \alpha\end{array}\right)\,\rangle_A=\int_M\,tr \,a\wedge \lambda,\label{canonical1}
\end{equation}
for any tangent vector \(\left(\begin{array}{c}a\\ \alpha\end{array}\right)\in T_{(A,\lambda)}S\). 

Let \(\phi\) be a 1-form on \({\cal A}\).   By definition,  \(\phi\) is a section of the cotangent bundle \(T^{\ast}{\cal A}\), so the pullback by \(\phi\) of \(\theta\) is a 1-form on \({\cal A}\).   We have the following characteristic property:
\begin{equation}
\phi^{\ast}\theta=\phi.
\end{equation}

\begin{lem}\label{lem1-1}
The derivation of the 1-form \(\theta\,\); 
 is given by
\begin{equation}
(\partial_{(A,\lambda)}\,\theta )\left(\begin{array}{c}a\\ \alpha\end{array}\right) = \int_M\,tr\,a\wedge \alpha
\,,\quad \mbox{ for }\,
\forall \left(\begin{array}{c}a\\ \alpha\end{array}\right)\in T_{(A,\lambda)}S\, .
\end{equation}
\end{lem}

In fact, 
\[
(\partial_{(A,\lambda)}\theta)(\left(\begin{array}{c}a\\ \alpha\end{array}\right))
=
\lim_{t\longrightarrow 0}\frac{1}{t}\int_M\,(\,tr\,a\wedge (\lambda+t\alpha)\,-\,tr\,a\wedge \lambda\,)
=
\int_M\,tr\,a\wedge \alpha
.
\]

The canonical 2-form is defind by 
\begin{equation}
\omega=\widetilde{d}\theta
.\end{equation}

 Lemma \ref{lem1-1}  and  (\ref{extder1}) yields the following
\begin{prop}
\begin{equation}
\omega_{(A,\lambda)}\,(\left(\begin{array}{c}a\\ \alpha\end{array}\right),\,\left(\begin{array}{c}b\\ \beta\end{array}\right)\,)=\int_M\,tr[ \,b\wedge\alpha-a\wedge\beta\,]
\end{equation}
\end{prop}
\(\omega\) is a non-degenerate closed 2-form on the cotangent space \(S\).     We see the non-degeneracy as follows.   Let \(\left(\begin{array}{c}a\\ \alpha\end{array}\right)\in T_{(A,\lambda)}T^{\ast}{\cal A}\), then \(a\in\Omega^1(M.Lie\,G)\) and \(\alpha\in \Omega^{m-1}(M,Lie\,G)\).   Hence 
\(\ast\alpha\in\Omega^1(M.Lie\,G)\) and \(\ast a\in\Omega^{m-1}(M.Lie\,G)\) and we have 
\[\omega_{(A,\lambda)}\,(\left(\begin{array}{c}a\\ \alpha\end{array}\right),\,\left(\begin{array}{c}\ast \alpha\\ \ast  a\end{array}\right)\,)=||\alpha||^2_{m-1}-||a||^2_1,\]
where \(\Vert\cdot\Vert_k\) is the \(L^2\)-metric on \(\Omega^k(M,\,Lie\,G)\).
This formula implies the non-degeneracy of \(\omega\).

Let \(\Phi=\Phi(A,\lambda) \) be a function on the cotangent space \(S\).   The Hamitonian vector field 
\(X_{\Phi}\) of \(\Phi\) is defined by the formula:
\begin{equation}
(\,\widetilde d\,\Phi\,)_{(A,\lambda)}\,=\,\omega (\quad\cdot\quad,\,X_{\Phi}(A,\lambda)\,).
\end{equation}

  Here the directional derivative of \(\Phi\) at the point \((A,\lambda)\) to the direction \(\left(\begin{array}{c}a\\\alpha\end{array}\right)\,\in T_{(A,\lambda)}S\) is  defined by the formula
\begin{equation}\label{partial2}
( \partial\Phi)_{(A,\lambda)}\left(\begin{array}{c}a\\\alpha\end{array}\right)\,=\lim_{t\longrightarrow 0}\,\frac{1}{t}(\Phi(A+ta,\lambda+t\alpha)-\Phi(A,\lambda)).
 \end{equation}
Hence the partial derivatives \((\frac{\delta\Phi}{\delta A})_A\in T^\ast_A\mathcal{A}\) and  \((\frac{\delta\Phi}{\delta \lambda})_A\in T_A\mathcal{A}\)  are given respectively by the formulas
\begin{eqnarray}
 \langle\,\frac{\delta\Phi}{\delta \lambda}\,,\,\alpha\, \rangle_{A}\,&=&\lim_{t\longrightarrow 0}\,\frac{1}{t}(\Phi(A,\lambda+t\alpha)-\Phi(A,\lambda)),\\[0.2cm]
 \langle\,a\,,\,\frac{\delta\Phi}{\delta A}\,\rangle_{A}\,&=&\lim_{t\longrightarrow 0}\,\frac{1}{t}(\Phi(A+ta,\lambda)-\Phi(A,\lambda)).\end{eqnarray}
It holds that 
\begin{equation}
(\,\widetilde d\,\Phi\,)_{(A,\lambda)}\left(\begin{array}{c}a\\ \alpha\end{array}\right)\,=\,
\langle\,a\,,\,\frac{\delta\Phi}{\delta A}\,\rangle_{A}\,+\, \langle\,\frac{\delta\Phi}{\delta \lambda}\,,\,\alpha\, \rangle_{A}\,.
\end{equation}
So the Hamiltonian vector field of \(\Phi\) is given by
\begin{equation}
X_{\Phi}\,=\,\left(\begin{array}{c}\, \,\frac{\delta\Phi}{\delta \lambda}\\[0.2cm]
-\, \frac{\delta\Phi}{\delta A}\end{array}\right).
\end{equation}

{\bf Example}

Let \(M\) be a compact three dimensional manifold.    We look at the following Hamiltonian function 
\begin{equation}
H(A,B)\,=\, \frac12\int_M\,tr[F_A\wedge\ast F_A\,]+\,\frac12\int_M\,tr[B\wedge \ast B]\,,
\end{equation}
for \(A\in \mathcal{A}\,\), \(B\in T^{\ast}_A\mathcal{A}\). 
   Then, since \(\frac{\delta H}{\delta A}=d_A\ast F_A=\ast ( d_A^{\ast}F_A)\,\) and \(\frac{\delta H}{\delta B}=\ast B\,\),  the Hamiltonian vector field of \(H\) becomes 
\begin{equation}
X_{H}\,=\,\left(\begin{array}{c}  \ast B \\[0.1cm]
- \ast\,  d_A^{\ast }F_A \end{array}\right).
\end{equation}
The Hamilton's equation of motion is  
\begin{equation}\label{*hamiltoneq}
\dot{A}\,=\,\ast B\,,\quad \dot{B}\,=\,- \ast \,d_A^\ast F_A\,, \qquad\mbox{ for } A\in \mathcal{A}\,,\,B\in T_A\mathcal{A}\,.
\end{equation}
It follows that the critical points of the Hamiltonian function \(H(A,0)= \frac12\int_Mtr[F_A\wedge\ast F_A\,]\)  are given by the Yang-Mills equation on the 3-dimensional manifold: 
\begin{equation}
d_AF_A= d_A^{\ast}F_A=0\,.
\end{equation}
\hfil\qed

The group of ( pointed ) gauge transformations \({\cal G}(M)=\Omega^0_s(M,\,Ad\,G)\) acts on \(T_A{\cal A}\)  by the adjoint representation; 
\(a\longrightarrow Ad_{g^{-1}}\,a=g^{-1}ag\), and on  \(T^{\ast}_A{\cal A}\) by its dual \(\alpha\longrightarrow g^{-1}\alpha g\).     Hence the canonical 1-form and 2-form are \({\cal G}\)-invariant.   
   The infinitesimal action of \(\xi\in Lie\,{\cal G}=\Omega^0(M,Lie\,G)\) on the cotangent space \(S=T^{\ast}{\cal A}\)  
gives a vector field \(\xi_S\) ( called fundamental vector field ) on \(S\) that is defined at the point \((A,\lambda)\) by the equation:
\begin{equation}
\xi_S(A,\lambda)=
\frac{d}{dt}\,\exp\,t\xi\cdot  \left(\begin{array}{c}A\\[0.2cm]
 \lambda\end{array}\right)=\left(\begin{array}{c}d_A\xi\\ [0.2cm]
\left[\lambda,\xi \right]\end{array}\right).
\end{equation}
Remember that  \({\cal G}_0(M)\) is the group of gauge transformations that are identity on the boundary of \(M\).       
When \(M\) has no boundary \({\cal G}(M)={\cal G}_0(M)\).  
If \(\xi\in  Lie\,{\cal G}_0\,\),  the vector field 
\(\,\xi_{T^{\ast}{\cal A}}(A,\lambda)\) is null on the boundary.

The dual space of \(\,Lie\,\mathcal{G}_0\) is given by \( (Lie\,\mathcal{G}_0)^{\ast}\,=\,\Omega^{m}(M, Lie\,G)\) with the dual pairing:
\begin{equation}
\langle \mu,\,\xi\rangle\,=\,\int_M\,tr\,(\mu\, \xi\,)\,,\qquad \forall\xi\in \mathcal{G}_0\,,\,\mu\in \Omega^m(M,Lie\,G).
\end{equation}

The {\it moment map} for the action of  \(\mathcal{G}_0\) on the symplectic space \((S\,,\omega)\) is the map 
\(\,
J: \,S\,\longrightarrow (Lie\,\mathcal{G}_0)^{\ast}\,\simeq \Omega^m(M,\,Lie\,G)\,,
\)
such that, if we denote 
\(J^{\xi}(A,\lambda)=\langle J(A,\lambda) ,\xi\rangle\,\) for \(\xi\in Lie\,\mathcal{G}_0\,\), 
\begin{enumerate}
\item
\(J^{\xi}\) is \(Ad^{\ast}\mathcal{G}\)-equivariant:
\begin{equation}
J^{Ad_g\xi}(g\cdot A,\,g\cdot\lambda)\,=\,J^{\xi}(A,\lambda)\,,
\end{equation}
\item
\(J^{\xi}\) satisfies the relation 
\begin{equation}\label{momentJ}
\widetilde{d}\,J^{\xi}\,=\,\omega\,(\,\,\cdot\quad,
\, \xi_S\,)\,.
\end{equation}
\end{enumerate}

\begin{prop}~~\\
The action of the group of gauge transformations \({\cal G}_0(M)\) on the symplectic space \((\,S\,,\omega\,)\) is an hamiltonian action with the moment map given by 
\begin{equation}
J(A,\lambda)=\,-\,d_A\,\lambda\,. \label{cotagmoment}
\end{equation} 
\end{prop}

{\it Proof}

The equivariance of \(J^{\xi}\) follows easily.    We shall verify the condition (\ref{momentJ}).
 Stokes' theorem yields 
\[J^{\xi}(A,\lambda)=
\langle J(A,\lambda),\xi\,\rangle =-\int_M\,tr\,(d_A\lambda)\,\xi\,=\,\int_M\,tr\,(\,d_A\xi\wedge\,\lambda)\,
  .\] 
 Since 
\[\lim_{t\longrightarrow 0}\frac{1}{t} \int_M tr\,\left( d_{A+ta}\xi\wedge (\lambda+t\alpha)-d_A\xi\wedge \lambda\right) =\,\int_M\,tr\,(a \wedge [\,\xi\,,\lambda\,]+\,d_A\xi\wedge \alpha\,.\]
we have 
\begin{equation}
\left(\widetilde{d}\,J^{\xi}\right)_{(A,\lambda)}\left(\begin{array}{c}a\\ \alpha\end{array}\right)\,=\,
\omega_{(A,\lambda)}\left(
\,\left(\begin{array}{c}a\\ \alpha\end{array}\right)\,, \left(\begin{array}{c}d_A\xi\\
\, [ \lambda\,, \xi ]\end{array}\right) \right)\,.
\end{equation}
\hfil\qed

\begin{defn}~~

The above moment map will be denoted by \(J^{\ast}:\,S\,\longrightarrow (Lie\,\mathcal{G})^{\ast}\,\).
\begin{equation}
J^{\ast}(A,\lambda)\,=\,- d_A\lambda\,,\qquad\forall (A,\lambda)\in S.
\end{equation}
\end{defn}

\begin{rem}~~

  \((J^{\ast})^{\xi}(A,\lambda)\) is given by 
 \begin{equation}
 \,(J^{\ast})^{\xi}(A,\lambda)\,=\,- \theta_{(A,\lambda)}\left(\xi_{T^{\ast}{\cal A}}\right)\,.
 \end{equation}
 \end{rem}
 \begin{rem}~~
From Marsden-Weinstein reduction theorem the reduced space \((J^\ast)^{-1}(\rho)/\mathcal{G}_0\) for a \(\rho\in (Lie\,\mathcal{G}_0)^{\ast}\) becomes a symplectic manifold endowed with the induced symplectic form \(\omega\),  and coincides with the space 
\(\{(A,\lambda)\in S\,;\,\, - d_A\lambda=\rho\,\}\).
 \end{rem}

\subsection{}

Let \(M\) be a {\it compact} \(m\)-dimensional riemannian manfold and \(G=SU(n)\), \(n\geq 2\), the special unirary group.   The Lie algebra of \(G\) is \(n\times n\)-matrices with vanishing trace.   Let \(P\,\stackrel{\pi}{\longrightarrow} M\) be the principal \(G\)-bundle over \(M\), and let \(\mathcal{A}\) be the space of irreducible connections over \(M\).  We assume that \(P\) is a trivial bundle \(P=M\times G\) though the same argument applies to the non-trivial case by a small change.    \(\mathcal{A}\) is an affine space modelled by the vector space \(\Omega^1(M,\,ad\,P\,=\Omega^1(M,\,su(n))\).   The tangent space at the point \(A\in\mathcal{A}\) is  
\begin{equation}
T_A\mathcal{A}\,=\,\Omega^1(M,su(n))\,.
\end{equation}
The inner product on \(T_A\mathcal{A}\) is given by 
\begin{equation}
(\,a\,,\,b\,)_1\,=\,\int_M\,Tr\,a\wedge\,\ast b\,\qquad \forall a,b\in T_A\mathcal{A}\,.\end{equation}
This is expressed by the product of differential forms and the multiplication of matrices.

We put \(R=T\mathcal{A}\).   The point of \(R\) is denoted by \((A,p)\) with  \(A\in\mathcal{A}\) and \(p\in T_A\mathcal{A}\).     The tangent space at \((A,p)\in R\) is \(T_{(A,p)}R=T_A\mathcal{A}\oplus T_A{\mathcal A}\), so any tangent vector \(\mathbf{a}\in T_{(A,p)}R\) is of the form \(\mathbf{a}=\left(\begin{array}{c}a\\x\end{array}\right)\,\) with \(a,\,x\in T_A\mathcal{A}\).   
The symplectic structure on \(R\) is defined by the formula 
\begin{equation}
\sigma_{(A,p)}\left(\,\left(\begin{array}{c}a\\x\end{array}\right)\,,\,\left(\begin{array}{c}b\\y\end{array}\right)\,\right)\,=\, (\,b\,,\,x\,)_1\,-\,(\,a\,,\,y\,)_1\,,\end{equation}
for all  \(\left(\begin{array}{c}a\\x\end{array}\right)\,,\,\left(\begin{array}{c}b\\y\end{array}\right)\,\in T_{(A,p)}R\,\).

The directional derivative to the direction \(\left(\begin{array}{c}a\\x\end{array}\right)\) of a function \(\varphi\) over \(R\) is defined by 
\begin{equation}
(\partial\,\varphi)_{(A,p)}\left(\begin{array}{c}a\\x\end{array}\right)\,=\,\lim_{t\longrightarrow 0}\,\frac{1}{t}\left(
\varphi(A+ta,p+ta)-\varphi(A,p)\right).
\end{equation}
Hence the partial derivatives \(\,\frac{\delta\varphi}{\delta A},\, \frac{\delta\varphi}{\delta p}\in T_A\mathcal{A}\simeq \Omega^1(M,su(n))\)
is defined by the equations:
\begin{equation}\label{partial1}
\partial\, \varphi_{(A,p)}\left(\begin{array}{c}a\\0\end{array}\right)\,=\,\left(\frac{\delta \varphi}{\delta A}\,,\, a\right)_1\,,
\qquad\partial\,  \varphi_{(A,p)}\left(\begin{array}{c}0\\x\end{array}\right)\,=\,\left(\frac{\delta \varphi}{\delta p}\,, \,x\right)_1\,.
\end{equation}
The hamiltonian vector  field \(X_\varphi\) corresponding to \(\varphi\) is given by 
\begin{equation}\label{Rhamiltonian}
X_\varphi\,=\,\left(\begin{array}{c}  \frac{\delta \varphi}{\delta p} \\[0.2cm] - \frac{\delta \varphi}{\delta A} \end{array}\right)\,.
\end{equation}

{\bf Example}
 
If we take the hamiltonian function
\begin{equation}
H(A,p)\,=\,\frac12(F_A,F_A)_2\,+\,\frac12(p,p)_1\,.
\end{equation}
Then
\begin{equation}
\partial H_{(A,p)}\left(\begin{array}{c}a\\x\end{array}\right)\,=\,(d_Aa,F_A)_2\,+\,(p,x)_1\,=\,(a,d_A^{\ast}F_A)_1+(p,x)_1\,,
\end{equation}
and 
 the corresponding hamiltonian vector field is 
\begin{equation}
(\,X_H\,)_{(A,p)}\,=
\,\left(\begin{array}{c}  p \\[0.2cm] - d_A^{\ast}F_A\end{array}\right)\,.
\end{equation}
Therefore the Hamilton's equation of motion becomes
\begin{equation}\label{hamiltoneq}
\dot A\,=\,p\,,\qquad 
\dot p\,=\,- d_A^{\ast}F_A \,.
\end{equation}

The group of gauge transformation \(\mathcal{G}=Aut_0(P)=\Omega^0(M,Ad P)\) acts on the symplectic manifold \((R,\omega)\):
\begin{equation}
g\cdot (A,p)\,=
\,(\,A+g^{-1}d_Ag\,,\,g^{-1}pg\,)\,, \quad g\in \mathcal{G}.
\end{equation}
Under the action of  \(\mathcal{G}\) the hamiltonian \(H\) is invariant.     

\begin{prop}\label{tgmoment}
The 
action of  \(\,\mathcal{G}\) on the symplectic space \((\,R=T\mathcal{A},\,\sigma\,)\) is an hamiltonian action with the moment map \(J:\,R\,\longrightarrow\,(Lie\, \mathcal{G})^{\ast}\simeq Lie\, \mathcal{G}\) given by
\begin{equation}
J(A,p)\,=\,d_A^{\ast}\,p\,. \label{tgmoment}
\end{equation}
\end{prop}

{\it Proof}~~

We note that we consider the dual of \(Lie\,\mathcal{G}\)  as  \(Lie\,\mathcal{G}=\Omega^0(M, Lie\,G)\) itself given by the inner product \((\,,\,)_0\).   
We shall prove that \(J^{\xi}(A,p)\,=\,(\,d_A^{\ast}p\,,\,\xi\,)_0\) gives the the moment map \(J(A,p): Lie\,\mathcal{G}\ni \xi\longrightarrow J^{\xi}(A,p)\in \mathbf{R}\) of the action of \(\mathcal{G}\) on \(R\).   
The fundamental vector field \(\xi_R\) corresponding to \(\xi\in Lie\,\mathcal{G}\,\) is 
\[\xi_R(A,p)\,=\,\frac{d}{dt}\vert_{t=0}(\,\exp\,t\xi\,\cdot A\,,\,\exp\,t\xi\cdot\,p\,)\,=\,\left(\begin{array}{c}d_A\xi
\\ -[\xi,\,p]\end{array}\right) \,\in \,T_{(A,p)}R\,.\]
Then, for any \(\left(\begin{array}{c}a\\x\end{array}\right)\in T_{(A,p)}R\),
\[\sigma(\,\left(\begin{array}{c}a\\x\end{array}\right)\,,\,\xi_R(A,p)\,)=\,(\,d_A\xi,x)_1\,-\,(a,-[\xi,p])_1\,.\]
On the other hand, since
\begin{eqnarray*}
(\widetilde d J^{\xi})_{(A,p)}\left(\begin{array}{c}a\\0\end{array}\right)&=&
\lim_{t\longrightarrow 0}\frac1t\left((d^\ast_{A+ta}p,\xi)_0-(d^\ast_Ap,\xi)_0\right)= \lim_{t\longrightarrow 0}\frac1t(\,p\,,\,d_{A+ta}\xi-d_A\xi)_1\\
&=&(p,[a,\xi])_1=(a,[\xi,p])_1\,,\end{eqnarray*}
and
\[(\widetilde dJ^\xi)_{(A,p)}\left(\begin{array}{c}0\\x\end{array}\right)=
\lim_{t\longrightarrow 0}\frac1t\left((d^\ast_{A}(p+tx),\xi)_0-(d^\ast_Ap,\xi)_0\right)
=(d_A^\ast x,\xi)_0=(x,d_A\xi)_1\,,\]
we have 
\[(\widetilde dJ^{\xi})_{(A,p)}\left(\begin{array}{c}a\\x\end{array}\right)
=(d_A\xi,x)_1+(a,[\xi,\, p])_1.\]
Therefore 
\[(\widetilde dJ^{\xi})_{(A,p)}\left(\begin{array}{c}a\\x\end{array}\right)\,=\,\sigma(\,\left(\begin{array}{c}a\\x\end{array}\right)\,,\,\xi_R(A,p)\,)\,,\quad \forall \left(\begin{array}{c}a\\x\end{array}\right)\in T_{(A,p)}R\,.\]
\(J(A,p)=d_A^{\ast}p\) is the moment map.
\hfil\qed

If \(\rho\in ( Lie\,\mathcal{G})^{\ast}\simeq \Omega^0(M,Lie\,G)\) then \(J^{-1}(\rho) =\{(A,p)\in R\,;\,d^{\ast}_Ap=\rho\,\}\).   So \(( Lie\,\mathcal{G})^{\ast}\) is the space of charge densities and the equation \(d^{\ast}_AE=\rho\,\) for   \(E=-p\) is the counterpart of Gauss's law.

\subsection{Duality}

Let \(M\) be a compact \(m\)-dimensional riemannian manifold.   Let \(\mathcal{A}\) be the space of irreducible connections over the trivial \(G=SU(n)\)-bundle \(P=M\times G\), and let \(\mathcal{G}\) be the 
group of gauge transformations over \(P\).
In previous sections we have investigated the symplectic manifolds  \((S=T^{\ast}\mathcal{A}\,,\,\omega)\)  and  \((R=T\mathcal{A}\,,\,\sigma)\).    There correspond Hamiltonian actions of the group of gauge transformations \(\mathcal{G}\) on these space with the moment maps \(J\), (\ref{tgmoment}), and \(J^{\ast}\), (\ref{cotagmoment}), respectively.    We must note that the dual space of \(Lie\,\mathcal{G}=\Omega^1(M,Lie\,G)\,\) viewed as the image of moment map \(J^{\ast}\) is \(\Omega^{m}(M,Lie\,G)\,\), while that as the image of \(J\) is \(\Omega^{1}(M,Lie\,G)\,\).

\begin{prop}   The symplectic manifolds \((S\,,\,\omega\,)\) and \((R\,,\,\sigma\,)\) are isomorphic via the Hodge operator:
 \begin{equation}
 \ast: T_A\mathcal{A}\,\simeq \Omega^1(M,su(n))\,\longrightarrow\,
T^{\ast}_A\mathcal{A}\,\simeq \Omega^{m-1}(M,su(n))\,, \quad A\in\mathcal{A}.\label{*iso}
\end{equation}
\end{prop}
In fact the correspondence  \( T_A\mathcal{A}\ni a\longrightarrow \ast a\in T^{\ast}_A\mathcal{A}\) gives the symplectic isomorphism:
\begin{equation*}
\sigma_{(A,p)}\left(\,\left(\begin{array}{c}a\\x\end{array}\right)\,,\,\left(\begin{array}{c}b\\y\end{array}\right)\,\right)= (\,b\,,\,x\,)_1\,-\,(\,a\,,\,y\,)_1\,=\,\omega_{(A,\ast p)}\left(\,\left(\begin{array}{c}a\\ \ast x\end{array}\right)\,,\,\left(\begin{array}{c}b\\ \ast y\end{array}\right)\,\right).\hfil\qed
\end{equation*}

\begin{prop}~~\label{orthgonaldecompo}
\begin{enumerate}
\item
\begin{enumerate}
\item
We have the following orthogonal decomposition of \(T_A\mathcal{A}\):
\begin{equation}
T_A\mathcal{A}\,=\,\{ d_A\xi\,;\,\xi\in Lie\,\mathcal{G}\,\}\oplus H_A\,,
\end{equation}
with \( H_A\,=\,\{x\in \Omega^1(M,\,Lie\,G);\, d^{\ast}_Ax=0\,\}\).
\item
Let 
\begin{equation}
R^0=\,\cup_{A\in \mathcal{A}}H_A\,.
\end{equation}  
Then \(R^0\) coincides with the symplectic reduction of \(R\) by the moment map \(J\),(\ref{tgmoment}):
\begin{equation}\label{reduction}
J^{-1}(0)/\mathcal{G}\,\simeq\, R^0\,.
\end{equation}
\end{enumerate}
\item 
\begin{enumerate}
\item
We have the following orthogonal decomposition of \(T^{\ast}_A\mathcal{A}\):
\begin{equation}
T^{\ast}_A\mathcal{A}\,=\,\{ d^{\ast}_A\lambda\,;\,\lambda\in ( Lie\,\mathcal{G})^{\ast}=\Omega^{m}(M,Lie\,G)\,\}\oplus H^{\ast}_A\,,
\end{equation}
with \( H^{\ast}_A\,=\,\{w\in \Omega^{m-1}(M,\,Lie\,G);\, d_Aw=0\,\}\).
\item
Let 
\begin{equation}
S^0=\,\cup_{A\in \mathcal{A}}H^{\ast}_A\,
\end{equation}
  Then 
 \(S^0\) coincides with the symplectic reduction of \(S\) by the moment map \(J^{\ast}\),(\ref{cotagmoment}):
\begin{equation}
 (J^{\ast})^{-1}(0)/\mathcal{G}\,\simeq\,S^0\,.
\end{equation}
\end{enumerate}
\item
The symplectic isomorphism given by Hodge \(\ast\), (\ref{*iso}), changes the 
Hamilton equations of motion (\ref{*hamiltoneq}) and (\ref{hamiltoneq})  together:
\begin{equation}
\{\begin{array}{c} \dot{A}\,=\,\ast B\\[0.2cm]  \dot{B}\,=\,- \ast\,d_A^\ast F_A\end{array}\quad
\stackrel{\,p=\ast B\,}\Longleftrightarrow\quad
\{\begin{array}{c} \dot A\,=\,p\\[0.2cm]
\dot p\,=\,- d_A^{\ast}F_A\end{array}\,.
\end{equation}
\end{enumerate}
\end{prop}

{\it Proof}

We shall prove only the orthogonal decomposition of \(T_A\mathcal{A}\) and that of \(T^{\ast}_A\mathcal{A}\).   The rest is easy to see.   Let \(H_A\) be the orthogonal complement of the subspace  \(d_A\,Lie\,\mathcal{G}\) in \(T_A\mathcal{A}=\Omega^1(M,Lie\,G)\).   
So  
\(x\in \Omega^1(M,\,Lie\,G)\) is decomposed to 
\(x=d_A\xi+y\,\) with \(\xi \in \Omega^0(M,\,Lie\,G)\) and \(y\in H_A\,\).   It implies \(H_A=\{y\in T_A\mathcal{A}:\, d_A^{\ast}y=0\}\).   It is easy to verify \( \xi=G_Ad^{\ast}_Ax\,\).    
Where \(G_A\) is the Green operator that will be precisely explained at the beginning of section 3.3.
   Next let \(u\in T^{\ast}{\cal A}=\Omega^{m-1}(M,\,Lie\,G)\) and put \(x=\ast u\).    Then 
 \(x\in \Omega^1(M,\,Lie\,G)\) is decomposed to 
\[x=d_A\xi+y\,,\quad\mbox{ with }\, \xi=G_Ad^{\ast}_Ax\,\in \Omega^0(M,\,Lie\,G),\quad y\in H_A\,.\]
Let \(\lambda=\ast\xi\,\), \(w=\ast y\,\).    Since \(d_Aw=\ast d_A^{\ast}y=0\) and \(d_A^{\ast}\lambda=
\ast d_A\xi\,\), we have \(u=d_A^{\ast}\lambda + w\) with \(w\in H^{\ast}_A\).   The decomposition is orthogonal with respect to \((\,\,,\,)_{m-1}\).
We have also \(\lambda=d_AG_Au\).

\hfil\qed

\section{Electronic Magnetic paradigm of Yang-Mills fields }

\subsection{Symplectic structure over \(T\mathcal{A}\times_{\mathcal{A}}T^{\ast}\mathcal{A}\)}

We introduce the  Yang-Mills field over the {\it 3-dimensional manifold} \(\,M\), that becomes the configuration space of the Yang-Mills equation.    This is a counterpart of the configuration space for the Hamiltonian description of Maxwell's equation, \cite{G}.   For that we consider the fiber product of the cotangent space \(T^{\ast}\mathcal{A}\) and the tangent space \(T\mathcal{A}\) over \(\mathcal{A}\).   In the following we shall endow it with a symplectic structure.

We introduce the Whitteney's direct sum of tangent and cotangent bundles: \(\,\mathbb{T}=T\mathcal{A}\times_{\mathcal{A}} T^{\ast}\mathcal{A}\longrightarrow\,\mathcal{A}\,\):  
\begin{equation}
\begin{array}{ccc}
\mathbb{T}=\,T\mathcal{A}\times_{\mathcal{A}}T^{\ast}\mathcal{A}\,&\,\longrightarrow\,&\,(T^{\ast}\mathcal{A}\,,\omega )=S\,\\[0.3cm]
\pi_{\ast}\,\downarrow\quad  &\qquad&\,\pi \downarrow\, \\[0.4cm]
R=(\,T\mathcal{A}\,,\sigma)\,&\,\longrightarrow\,&\, \,\mathcal{A}\,.
\end{array}
\end{equation}      
A point of \(\,\mathbb{T}\,\) will be denoted by \((A, E,B)\) with \(A\in \mathcal{A}\), \(E\in T_A\mathcal{A}\) and  \(B\in T^{\ast}_A\mathcal{A}\).   
The tangent space of \(\,\mathbb{T}\,\) at the point \((A,E,B)\in \mathbb{T}\,\)   is  \( T_{(A,E,B)}\mathbb{T}\,=T_A\mathcal{A}\oplus (T_A\mathcal{A}\oplus T^{\ast}_A\mathcal{A})\).   So a vector in \(T_{(A,E,B)}\mathbb{T}\,\) is denoted by \(\mathbf{a}=\left(\begin{array}{c}a\\ e\\ \beta\end{array}\right) \), with \(a,\,e\in T_A\mathcal{A}\simeq \Omega^1(M,Lie\,G)\) and \(\beta\in T_A^{\ast}\mathcal{A}\simeq \Omega^{2}(M,Lie\,G)\).        
Remember that each fibers \(T_A\mathcal{A}\simeq  \Omega^1(M,Lie\,G)\) and \(T_A^{\ast}\mathcal{A}\simeq \Omega^{2}(M,Lie\,G)\)  are endowed with the structure of Sobolev space \(L^2_{s-1}\) with \(s>2\).            

Now we shall endow the space \(\mathbb{T}\) with a symplectic structure that is written only by the field variables \((E,B)\), so the potential variable \(A\) is implicit.    

Let \(\Theta\) be the following 1-form on \(\mathbb{T}\) defined by
\begin{equation}\label{C1form}
\Theta_{(A,E,B)}\left(\left(\begin{array}{c}a\\e\\ \beta\end{array}\right)\right)\,=\,\left(\,e\,,\,d_A^{\ast}B\,\right)_1\,,
\quad\forall \left(\begin{array}{c}a\\e\\ \beta\end{array}\right)\in T_{(E,B)}\mathbb{T}.
\end{equation}

Since \((\partial_A\,d_A^{\ast}\,)_{(A,E,B)}a\,=\ast(a\ast\cdot)\), the partial derivative with respect to \(A\) of the 1-form \(\Theta\) to the direction \(a\in T_A\mathcal{A}\) is
\begin{equation}
(\partial_A\,\Theta)_{(A,E,B)}\left(\begin{array}{c}a\\0\\0\end{array}\right)\,\,=\,(\,\cdot\,,\,\ast (a\ast B)\,)_1\,=\,(\,\cdot\wedge a\,,\,B\,)_2\,,
\end{equation}
The partial derivative with respect to \(E\) of the 1-form \(\Theta\) to any direction vanishes;
\begin{equation}
(\partial_E\,\Theta)_{(E,B)}\,=\,0\,,
\end{equation}
while the partial derivative with respect to \(B\) of  \(\Theta\) is given by 
\begin{equation}
(\partial_B\,\Theta)_{(E,B)}\left(\begin{array}{c}0\\0\\ \beta\end{array}\right)\,=\,(\,\cdot\,,d_A^{\ast}\beta\,)_1.
\end{equation}
Hence the exterior derivative of the \(1\)-form \(\Theta\) becomes
\begin{eqnarray}
&&(\widetilde d\Theta)_{(E,B)}\left(\left(\begin{array}{c}a_1\\e_1\\ \beta_1\end{array}\right)\,,\,\left(\begin{array}{c}a_2\\e_2\\ \beta_2\end{array}\right)\,\right)\\[0.2cm]
\qquad&&=\,(e_2\,a_1\,,\,B\,)_2\,-\,(e_1a_2\,,\,B\,)_2\,+
\,(e_2,d_A^{\ast}\beta_1)_1\,-\,(e_1,d_A^{\ast}\beta_2)_1\,,
\end{eqnarray}
see the formula (\ref{extder1}) in about the exterior differentiation of a 1-form. 

\begin{defn}
We define the  2-form \(\Omega\) on  \(\mathbb{T}\) by the following formula:
\begin{equation}
\Omega\,=\,\widetilde d\,\Theta\,.
\end{equation}
\end{defn}
Then
\begin{equation}
 \Omega\,=\,\Omega^1\,+\,\Omega^2\,
\end{equation}
with
\begin{eqnarray}\label{vorsymp}
  \Omega^1_{(A,E,B)}\left(\left(\begin{array}{c}a_1\\e_1\\ \beta_1\end{array}\right)\,,\,\left(\begin{array}{c}a_2\\e_2\\ \beta_2\end{array}\right)\,\right)\,
&=&
\,(e_2\,a_1\,,\,B)\,-\,(e_1a_2\,,\,B)\,,\nonumber\\[0.2cm]
 \Omega^2_{(A,E,B)}\left(\left(\begin{array}{c}a_1\\e_1\\ \beta_1\end{array}\right)\,,\,\left(\begin{array}{c}a_2\\e_2\\ \beta_2\end{array}\right)\,\right)\,
& =&
\,(e_2\,a_1\,,\,B)\,-\,(e_1a_2\,,\,B)\,\nonumber
\end{eqnarray}
Since
 \[(\partial_B\Omega^1)_{(A,E,B)}\beta_1\left(\left(\begin{array}{c}a_2\\e_2\\ \beta_2\end{array}\right)\,,\,\left(\begin{array}{c}a_3\\e_3\\ \beta_3\end{array}\right)\,\right)\,
 =\left(e_3a_2\,,\,\beta_1\right)\,-\,\left(e_2a_3\,,\,\beta_1\,\right),\]
 and
 \[(\partial_A\Omega^2)_{(A,E,B)}\,a_1\left(\left(\begin{array}{c}a_2\\e_2\\ \beta_2\end{array}\right)\,,\,\left(\begin{array}{c}a_3\\e_3\\ \beta_3\end{array}\right)\,\right)\,
 =\left(e_3a_1\,,\,\beta_2\right)\,-\,\left(e_2a_1\,,\,\beta_3\,\right),\]
 and other directional derivations are \(0\), we have 
\(\widetilde d\,\Omega=\widetilde d\,\Omega^1+\widetilde d\,\Omega^2=0\) from (\ref{extder}).     It is easy to see that \(\Omega\) is a non-degenerate skew linear form.   Thus we have the following
\begin{thm}~~\\
\(\left(\, \mathbb{T}=T\mathcal{A}\times_{\mathcal{A}}T^{\ast}\mathcal{A}\,,\,\Omega\,\right)\) is a symplectic manifold.
\end{thm}

\vspace{0.2cm}

We define the inner product  on each fiber of \(\,\mathbb{T}\longrightarrow\mathcal{A}\).         
\begin{equation}
\left( \left(\begin{array}{c}e_1\\ \beta_1\end{array}\right) ,\,\left(\begin{array}{c}e_2\\ \beta_2\end{array}\right) \right)_{\mathbb{T}}\,=\, (e_2,d_A^{\ast}\beta_1\,)_1\,+\,(e_1,\,d_A^{\ast}\beta_2\,)_1\,,
\end{equation}
for \(E\in T_A\mathcal{A}\), \(B\in T_A^{\ast}\mathcal{A}\) and \(\forall  \left(\begin{array}{c}e_i\\ \beta_i\end{array}\right) \in T_{(A,E,B)}\mathbb{T}\,\), \(i=1,2\), and 
.    Here a tangent vector to the fiber \(T_A\mathcal{A}\) is denoted by \( \left(\begin{array}{c}e_1\\ \beta_1\end{array}\right)\).

The partial derivaive of a function \(\Phi=\Phi(E,B)\) over \(\mathbb{T}\) to the direction \(e\in T\mathcal{A}\) 
 is defined as  the vector  \(\frac{\delta\Phi}{\delta E}\,\in T_A^{\ast}\mathcal{A}\,\) that satisfies the equation
 \begin{equation}
 (\partial\Phi)_{(E,B)}\left(\begin{array}{c}e\\0\end{array}\right)\,=\,(e,d_A^{\ast}\frac{\delta\Phi}{\delta E}\,)_1\,,\quad \forall e\in T_A\mathcal{A},
 \end{equation}
where \(\pi_{\ast}(E,B)=A\).
Respectively the directional derivative  to  \(\beta\in T^{\ast}\mathcal{A}\) is defined as  the vector  \(\frac{\delta\Phi}{\delta B}\,\in T_A\mathcal{A}\,\) that satisfies the equation
 \begin{equation}
 (\partial\Phi)_{(E,B)}\left(\begin{array}{c}0\\\beta\end{array}\right)\,=\,(\frac{\delta\Phi}{\delta B}\,,d_A^{\ast}\beta\,)_1\,,\quad \forall \beta \in T^{\ast}_A\mathcal{A}.
 \end{equation}
 We have then 
 \begin{equation}\label{dPhi}
( \widetilde d\Phi)_{(E,B)}\left(\begin{array}{c}e\\\beta\end{array}\right)\,=\,
\left( \left(\begin{array}{c}\frac{\delta\Phi}{\delta B}\\[0.2cm] \frac{\delta\Phi}{\delta E}\end{array}\right) ,\,\left(\begin{array}{c}e\\[0.2cm] \beta\end{array}\right) \right)_{\mathbb{T}}\,
\end{equation}
 
 {\bf Example}~
 
 Let \(H=H(E,B)\) be the Hamiltonian function on \(\mathbb{T}\) written in the vortex representation;
 \begin{equation}\label{vortexHamiltonian}
 H(E,B)\,=\frac12 \left(
 \left(\begin{array}{c}E\\ B\end{array}\right)\,,
 \left(\begin{array}{c}d_A^{\ast}B\\ d_AE\end{array}\right)
 \right)_{\mathbb{T}} =\frac12\{(d_AE,d_AE)_1+(d_A^{\ast}B,d_A^{\ast}B)_1\},
 \end{equation}
 we have
 \begin{equation}\label{Maxwelleq}
 \frac{\delta H}{\delta B}=d_A^{\ast}B\,,\quad  \frac{\delta H}{\delta E}=d_AE\,.
 \end{equation}
 \hfil\qed

\begin{prop}~~\\
Let \(\Phi=\Phi(E,B)\) be a function on the fields \(\mathbb{T}\).   
Then the Hamiltonian vector field \(X_{\Phi}\) of \(\Phi\) is given by 
\begin{equation}\label{fieldHam}
X_{\Phi}(E,B)\,=\,\left(\begin{array}{c}- \frac{\delta\Phi}{\delta B}\\[0.2cm]
\frac{\delta\Phi}{\delta E}\end{array}\right)\,.
\end{equation}
\end{prop}
The formulae  (\ref{dPhi}) and (\ref{vorsymp}) imply (\ref{fieldHam}).

\begin{defn}~~\\
The Poisson bracket on \(\mathbb{T}\) is defined by the formula
\begin{equation}\label{fieldP}
\left\{\,\Phi\,,\,\Psi\,\right\}_{(E,B)}^{\mathbb{T}}\,=\,\Omega_{(E,B)}\left(\,X_{\Phi}\,,\,X_{\Psi}\,\right)\,,
\end{equation}
for \(\Phi,\,\Psi\in\,C^{\infty}(\mathbb{T})\).
\end{defn}

The following formula is our counterpart to the Marsden-Weinsrein's vortex formula for the Poisson baracket of Maxwell's fields,\cite{MW2} .

\begin{prop}
\begin{equation}\label{MW}
\begin{aligned}
\left\{\,\Phi\,,\,\Psi\,\right\}_{(E,B)}^{\mathbb{T}}\,&=\,\left(\,\frac{\delta\Phi}{\delta B}\,,\,d_A^{\ast}\frac{\delta\Psi}{\delta E}\,\right)_1\,-
\,\left(\,\frac{\delta\Psi}{\delta B}\,,\,d_A^{\ast}\frac{\delta\Phi}{\delta E}\,\right)_1\\[0.3cm]
\quad&=\,\left(\,d_A\frac{\delta\Phi}{\delta B}\,,\,\frac{\delta\Psi}{\delta E}\,\right)_2\,-
\,\left(\,d_A\frac{\delta\Psi}{\delta B}\,,\,\frac{\delta\Phi}{\delta E}\,\right)_2
\end{aligned}
\end{equation}
\end{prop}

Proposition follows from (\ref{fieldHam}).

{\bf Example}

Let \(H=H(E,B)\) be the Hamiltonian function of (\ref{vortexHamiltonian}).    
The equation of motion on the strength field \(\mathbb{T}\) is written in the form 
\[\dot{\Phi}\,=\, \left\{\,\Phi\,,\,H\right\}_{(E,B)}^{\mathbb{T}}\,.\]
Which implies by the formula  (\ref{Maxwelleq})  the Maxwell equation over \(\mathbb{T}\):
\begin{equation}\label{YMMaxwell}
\dot E\,=\,-d_A^{\ast}B\,,\qquad \dot B\,=\,d_AE\,.
\end{equation}
\hfil\qed

The group of gauge transformations \(\mathcal{G}\) acts on \(\mathbb{T}\) by 
\begin{equation}
g\cdot (E,B)\,=\, (Ad_g\,E\,,\,Ad_{g^{-1}}^{\ast}B\,)\,=\,(g^{-1}Eg\,,\,g^{-1}Bg\,)\,.
\end{equation}
It is a symplectic action because of 
\[\left (\,g\cdot e\,,\,d_{g\cdot A}^{\ast}(g\cdot \beta)\, \right)_1\,=\, \left (\,g\cdot e\,,\,g\cdot (d_{ A}^{\ast} \beta \,\right)_1\,=\,\left( e\,,\,d_A^{\ast}\beta\,\right),\]
for any \((e,\beta)\in T_{(E,B)}\mathbb{T}\).     

The Lie algebra of infinitesimal gauge transformations  \(Lie\,\mathcal{G}=\Omega^0(M,Lie\,G)\) acts on \(\mathbb{T}\) by 
\begin{equation}   
\xi\cdot \left(\begin{array}{c}A\\[0.2cm]E\\[0.2cm]B\end{array}\right)\,=\, \left(\begin{array}{c}d_A\xi\\[0.2cm][\,E,\xi\,]\\[0.2cm][\,B,\xi\,]\end{array}\right)\,,
\end{equation}\label{fundamentalvf}
that is, the fundamental vector field on \(\mathbb{T}\) corresponding to \(\xi\in Lie\,\mathcal{G}\) becomes 
\begin{equation}
\xi_{\mathbb{T}}\,(A,E,B)\,=\,\left(\begin{array}{c}d_A\xi\\[0.2cm] [\,E,\xi\,]\\[0.2cm][\,B,\xi\,]\end{array}\right)\,.
\end{equation}

\subsection{Yang-Mills fields}

We shall introduce the Yang-Mills field as a subspace of \(\mathbb{T}\).      We have investigated the moment map over the space of vector potentials \(T\mathcal{A}\) and that over the space of vortex potentials \(T^{\ast}\mathcal{A}\):
\[J:\,R=T\mathcal{A}\,\longrightarrow (Lie\,\mathcal{G})^{\ast}\simeq \Omega^0(M, Lie\,G),\]
and
\[J^{\ast}:\,T^{\ast}\mathcal{A}\longrightarrow (Lie\,\mathcal{G})^{\ast}\,\simeq \Omega^3(M, Lie\,G).\]
 For a Hamiltonian system with symmetries the range of a moment map represents conserved quantities of the system.    In our case, we take \((Lie\, \mathcal{G})^{\ast}\simeq\Omega^3(M,Lie\,G)\) as the space of charges for the action \(\mathcal{G}\) over \(\mathbb{T}\),  that has two component corresponding to the electric charge ( times the volume form \(\ast 1\)) and the magnetic charge.   So given  \(\rho\ast 1\,,\,\rho^{\prime}\in (Lie\,\mathcal{G})^{\ast}\),  we have the configuration space  \(\{(E,B)\in\mathbb{T};\quad J(E)=\rho,\,J^{\ast}(B)=\rho^{\prime}\}\).     On the component which contains the field stength \(F_A\) 
we must have \(\rho^{\prime}=0\).

\begin{defn}~~

The {\it Yang-Mills field} is  the subspace of \(\mathbb{T}\) defined by 
\begin{equation}\label{field}
\mathcal{F}\,=\,\left\{(E,B)\,\in  \mathbb{T}\,:\quad d_AB=0,\,\,d_A^{\ast}E=0 \quad\mbox{ with \(\,\pi_{\ast}(E,B)=A\)}\,\right\}\,
\end{equation}
\end{defn}
The Yang-Mills field \(\mathcal{F}\) is a symplectic subspace of \((\mathbb{T},\Omega)\) that is \(\mathcal{G}\)-invariant  because of the relation:
\[d_{g\cdot A}(g\cdot B)\,=\,g\cdot(d_AB)\,,\quad d_{g\cdot A}(g\cdot E)\,=\,g\cdot(d_AE)\,.\]

On Yang-Mills field \(\mathcal{F}\) we have the analogous formula of Maxwell's equations:
\begin{eqnarray}\label{YMaxwell}
d_A^{\ast}B\,+\,\dot E\,=\,0\,&,&\quad d_A^{\ast}E\,=0\,,\\[0.2cm]
d_AE\,-\,\dot B\,=\,0\,&,&\quad d_AB\,=\,0\,.
\end{eqnarray}

Now we shall investigate the Hamiltonian action of \(\mathcal{G}\) on the Yang-Mills field \(\mathcal{F}\).

A map 
\begin{equation}
\mathbb{J}\,:\,\mathcal{F}\,\longrightarrow\,(Lie\,\mathcal{G})^{\ast}=\Omega^3(M,\,Lie\,G)
\end{equation}
is by definition a moment map for  the symplectic action of  \(\,\mathcal{G}\) on \(\,\mathcal{F}\,\) provided  
\begin{enumerate}
\item
If we put  \(\mathbb{J}^{\xi}(E,B)=\langle 
\mathbb{J}(E,B)\,,\xi\rangle\),
 the Hamiltonian vector fields of 
 \(\,\mathbb{J}^{\xi}\) coincides with the fundamental vector field \(\xi_{\mathbb{T}}\), (\ref{fundamentalvf}).
 \item
\(\mathbb{J}\) is \(Ad^{\ast}\)-equivariant:
\[\mathbb{J}^{\xi}(g^{-1}Eg\,,\,g^{-1}Bg\,)\,=\,\mathbb{J}^{Ad_{g^{-1}}\xi}(E,B)\,.\]
\end{enumerate}
In this case we say that the action of \(\mathcal{G}\) is Hamiltonian.

\begin{prop}
The action of  \(\,\mathcal{G}\) on \(\,\mathcal{F}\) is Hamiltonian with the moment map
\begin{equation}
\mathbb{J}(\,E,\,B\,)\,=\,[\,d_A\ast B\,,\,E\,]\,.
\end{equation}
\end{prop}

{\it Proof}

It holds the relation
\begin{equation}
\mathbb{J}^{\xi}(E,B)\,= \,\Theta_{(A,E,B)}(\,\xi_{\mathbb{T}}\,(A,E,B)\,).
\end{equation}

We have 
\begin{eqnarray*}
(\tilde d\,{\mathbb{J}}^{\xi})_{(A,E,B)}\,\left(
\begin{array}{c}a\\[0.2cm] e\\[0.2cm] \beta\end{array}\right)
&=&
([\,E,\xi]\wedge a\,,\,B)_2+([e,\xi]\,,\,d_A^{\ast}B\,)_1+(\,[E,\xi] \,,\,d_A^{\ast}\beta\,)_1\,,
\\[0.2cm]
\left(\,i_{\xi_{\mathbb{T}}}\,\Omega\,\right)\,\left(
\begin{array}{c}a\\[0.2cm] e\\[0.2cm] \beta\end{array}\right)\,&=&\, (\,[e,d_A\xi]\,,\,B)_2-(\,[E,\xi]\wedge a \,,\,B\,)_2\\[0.2cm]
\qquad&\,&\,+\,(e,d_A^{\ast}[B,\xi]\,)_1-(\,[E,\xi] \,,\,d_A^{\ast}\beta\,)_1\,.
\end{eqnarray*}
Hence \(\,\tilde d\,{\mathbb{J}}^{\xi}\,=\,\,-i_{\xi_{\mathbb{T}}}\,\Omega\,\).   
The equivariance of \(\mathbb{J}\) is easy to verify.
\hfil\qed

\begin{cor}
On the \(\mathcal{G}\)-orbit passing through a solution of  equation (\ref{YMaxwell}) 
we have \(\mathbb{J}(E,B)=\,[\dot B,\ast B]\).
\end{cor}

\subsection{Symplectic variable \(\gamma: R\longrightarrow \mathcal{F}\)}

Since any \(A\in \mathcal{A}\) is an irreducible connection we have the Green operator \(G_A\) defined on  \(\Omega^k((M,Lie\,G)\), \(k=1,2\), 
\[ (\,d_Ad^{\ast}_A+d^{\ast}_Ad_A\,)\,G_A\alpha\,=\,\alpha\,,\quad  \forall \alpha\in \Omega^k((M,Lie\,G).\]
  \(G_A\)  is a self adjoint operator; \((G_Au,v)_k=(u,G_Av)_k\) for any 
\(u,\,v\in \Omega^k(M, Lie\,G)\), \(k=1,2\).     We note also the fact that \(G_A\) commutes with \(d_A\) and \(d_A^{\ast}\,\):
\[ d_AG_A=G_Ad_A\,,\qquad   d^{\ast}_AG_A=G_Ad^{\ast}_A\,.\]

Restricted to the space \(\mathcal{F}\) we have 
\begin{equation}\label{Green}
d_Ad_A^{\ast}G_A\beta=\beta,\qquad d^{\ast}_Ad_AG_Ae=e,
\end{equation}
for \(e\in T_A\mathcal{A}\) and \(\beta\in T^{\ast}_A\mathcal{A}\).

\begin{defn}~~
\begin{enumerate}
\item
 \(\phi:\,R\,\longrightarrow\,\mathcal{F}\subset \mathbb{T}\,\)  is the map defined by 
\begin{equation}
\phi\left(A\,,\,p\right)\,=\,\left(E=-p\,,\, B=F_A\,\right).
\end{equation}
\item 
Let \(\phi_{\ast}:TR\longrightarrow T\mathbb{T}\) be the tangent map of \(\phi\,\):
\[(\phi_{\ast})_{(A,p)}\,\left(\begin{array}{c}a\\ x\end{array}\right)\,=\,\left(\begin{array}{c}-x\\ d_Aa\end{array}\right),\]
and let \(G_A:\,T_A\cal{A}\longrightarrow T_A\cal{A}\) be the Green operator.  \\
We define the modified tangent map   \(\gamma:\,TR\longrightarrow\,T\mathbb{T}\) of \(\phi\) as follows 
\begin{equation}
\gamma=\phi_{\ast} \circ \left(\begin{array}{cc}1&0\\[0.2cm]0& G_A\end{array}\right)=\left(\begin{array}{cc}0&-G_A\\[0.2cm]d_A& 0\end{array}\right)
,
\end{equation}
that is,
\[T_{(A,p)}R\,\ni\,\left(\begin{array}{c}a\\ x\end{array}\right)\,\longrightarrow\, \gamma_{(A,p)}\left(\begin{array}{c}a\\ x\end{array}\right)=
\left(\begin{array}{c}-G_Ax\\ d_Aa\end{array}\right)\in T_{\phi(A,p)}\mathbb{T}\,.\]
\end{enumerate}
\end{defn}

\begin{lem}\label{gamma-1}
\begin{equation}
\gamma^{\ast}\Omega\,=\,\sigma\,.
\end{equation}
\end{lem}
In fact, we have, for any \(\left(\begin{array}{c}a_i\\x_i\end{array}\right)\,\in T_{(A,p)}R\), \(i=1,2\),
\begin{eqnarray*}
(\gamma^{\ast}\Omega)_{(A,p)}
\left(
\left(\begin{array}{c}a_1\\x_1\end{array}\right)\,,
\,\left(\begin{array}{c}a_2\\x_2\end{array}\right)
\right)
&=&
\Omega_{(E,B)}
\left(
\gamma\left(\begin{array}{c}a_1\\x_1\end{array}\right)\,,
\gamma\left(\begin{array}{c}a_2\\x_2\end{array}\right)\,
\right)
=\\[0.3cm]
\Omega_{(E,B)}
\left(
\left(\begin{array}{c}-G_Ax_1\\d_Aa_1\end{array}\right)\,,
\left(\begin{array}{c}-G_Ax_2\\d_Aa_2\end{array}\right)\,
\right)
&=&
(-G_Ax_2\,,\,d^{\ast}_Ad_Aa_1\,)_1\,-\,(-G_Ax_1\,,\,d^{\ast}_Ad_Aa_2\,)_1\\[0.3cm]
\,=\,(x_1\,,\,G_Ad^{\ast}_Ad_Aa_2\,)_1\,-\, (x_2\,,\,G_Ad^{\ast}_Ad_Aa_1\,)_1\,
&=&
\,(x_1,a_2)_1\,-\,(x_2,a_1)\\[0.3cm]
=\,\sigma_{(A,p)}
\left(
\left(\begin{array}{c}a_1\\x_1\end{array}\right)\,,\,\left(\begin{array}{c}a_2\\x_2\end{array}\right)
\right)\,.&\,&\,
\end{eqnarray*}
\hfil\qed

Let \(R^0=\cup_{A\in\mathcal{A}}\,H_A^0\) be the reduction of \(R\,\), (\ref{reduction}).      Remenber that the symplectic reduction of \(R\) by the moment map \(J\) is isomorphic to \(R^0\),  Proposition \ref{orthgonaldecompo} .   
 
  \begin{thm}~~
 \((R^0,\sigma)\) is symplectomorph to  \((\mathcal{F},\Omega)\).
 \end{thm}
 
 {\it Proof}
 
Since  \(d^{\ast}_A(-p)=0\) and \(d_AF_A=0\) for \((A,p)\in R^0\),  \(\phi\) maps the subspace \(R^0\) into \(\mathcal{F}\) .     The tangent space of \(H^0_A\)  consist of  those vectors 
\( \left(\begin{array}{c}a\\x\end{array}\right)\in T_{(A,p)}R\) such that \(\,d_A^{\ast}x=0\), and the tangent space  \(T\mathcal{F}\) consists of those vectors \( \left(\begin{array}{c} e \\ \beta\end{array}\right)\in T\mathbb{T}\) such that \(d^{\ast}_Ae=0\) and \(d_A\beta=0\).   
 If  \( \left(\begin{array}{c}a\\x\end{array}\right)\) is  tangent to \( H^0_A \)  then \(d_A^{\ast}G_Ax=0\) and \(d_A(d_Aa)=0\), ( the latter follows from the derivation of \(d_AF_A=0\)).    So   \(\gamma\) maps \(TR^0\) into \(T\mathcal{F}\). 
Moreover \(\gamma\) is  a bijective map of  \(TR^0\) onto \(T\mathcal{F}\) .   In fact we have the inverse map 
 given by 
\[ T\mathcal{F}\ni \left(\begin{array}{c}e\\ \beta\end{array}\right)\,\longrightarrow \,(\,- d_A^{\ast}\circ \gamma)\left(\begin{array}{c}e\\ \beta\end{array}\right)=\left(\begin{array}{c}d_A^{\ast}G_A\beta\\ -d_A^{\ast}d_A e\end{array}\right)\in TR^0
.\]
By virtue of the implicit function theorem in Banach space  the vector spaces \(R\) and \(\mathcal{F}\) are diffeomorphic.   Let    \(\tilde\gamma: \,R^0\longrightarrow\mathcal{F}\) be the diffeomorphism.   Lemma \ref{gamma-1} implies that \(\tilde\gamma\) is a symplectomorphism.

\hfil\qed

\begin{lem}
Let \(\Phi \in C^{\infty}(\mathcal{F})\) and let \(\varphi \in C^{\infty}(R)\) be the pullback of \(\Phi\)  by \(\phi:R\longrightarrow\mathcal{F}\): 
\[\varphi(A,p)=\phi^{\ast}\Phi=\Phi(-p,F_A).\]
Then the Hamiltonian vector field \(X^{R}_{\varphi}\,\) of \(\varphi\) has the formula
\begin{equation}\label{gamaHamiltonian}
X^{R}_{\varphi}(A,p)\,=\,\left(\begin{array}{c}d_A^{\ast}\frac{\delta\Phi}{\delta E}\\[0.2cm]
d_A^{\ast}d_A \frac{\delta\Phi}{\delta B}\end{array}\right)_{(E,B)=\phi(A,p)}\,
\end{equation}
\end{lem}

{\it Proof}~~\\
\begin{eqnarray*}
(\partial \varphi)_{(A,p)}\left( 
\begin{array}{c}a\\0\end{array} \right)
&=& \lim_{t\longrightarrow 0}\,\frac{1}{t}\left( \Phi(- p,F_{A+ta})-\Phi(- p,F_A) \right)\\[0.2cm]
&=&  \lim_{t\longrightarrow 0}\,\frac{1}{t} \left(
\Phi(E,B+t\beta)\,-\,\Phi(E,B)
\right)\lvert_{
E=- p,\,B=F_A,\,\beta=d_Aa
}
\\[0.2cm]
&=&(\partial\Phi)_{(- p,F_A)}\left(\begin{array}{c}0\\ d_Aa\end{array}\right)
=
\,\left(\,\frac{\delta\Phi}{\delta B}\,,\,d_A^{\ast}d_Aa\,\right)_1\vert_{(-p,F_A)}\,\\[0.2cm]
&=&\,\left( a\,,\,(\,d_A^{\ast}d_A\frac{\delta\Phi}{\delta B}\,)\vert_{(- p,F_A)}\,\,\right)_1\,.
\end{eqnarray*}
Hence \(\frac{\delta\varphi}{\delta A}=\,(d_A^{\ast}d_A\frac{\delta\Phi}{\delta B}\,)\vert_{(- p,F_A)}\).
And 
\begin{eqnarray*}
(\partial \varphi)_{(A,p)}\left( 
\begin{array}{c}0\\x\end{array} \right)
&=&
\lim_{t\longrightarrow 0}\,\frac{1}{t}\left( \Phi(E+te),B)-\Phi(E,B) \right)
\vert_{E=-p,\,e=-x,\,B=F_A\,}\\[0.2cm]
&=&\left(\,-x\,,
\,(\,d_A^{\ast}\frac{\delta\Phi}{\delta E}\,)\vert_{(-p,F_A)}\,\right)_1\,=
\left(\,-(\,d_A^{\ast}\frac{\delta\Phi}{\delta E})\vert_{(-p,F_A)}\,,\,x\right)_1\,.
\end{eqnarray*}

Therefore
\[
\frac{\delta\varphi}{\delta A}\,=\,(\,d_A^{\ast}d_A\frac{\delta\Phi}{\delta B}\,)\vert_{(-p,F_A)},\quad
\frac{\delta\varphi}{\delta p}\,=\,-(\,d_A^{\ast}\frac{\delta\Phi}{\delta E}\,)\vert_{(-p,F_A)}\,.\]
(\ref{gamaHamiltonian}) follows from 
(\ref{Rhamiltonian}).
\hfil\qed

\begin{cor}
\begin{equation}
\gamma^{-1}X^{\mathbb{T}}_{\Phi}\,=\,
\,\left(\begin{array}{c}d_A^{\ast}G_A\frac{\delta\Phi}{\delta E}\\[0.2cm]
d_A^{\ast}d_A \frac{\delta\Phi}{\delta B}\end{array}\right)_{(E,B)=\phi(A,p)}\,
=\left(\begin{array}{cc}G_A&0\\[0.2cm]0&1\end{array}\right)\,X^{R}_{\varphi}.
\end{equation}
\end{cor}

\begin{prop}~~
\begin{equation}
\left\{\,\Phi,\,\Psi\,\right\}^{\mathbb{T}}\circ\tilde\gamma\,
=\,\sigma(
\,\gamma^{-1} X_{\Phi}^{\mathbb{T}},\,\gamma^{-1}X^{\mathbb{T}}_{\Psi}
\,)\,=\,\left\{\Phi\circ\tilde\gamma\,,\,\Psi\circ\tilde\gamma\,\right\}\,.
\end{equation}
\end{prop}

  \end{document}